\documentclass[12pt,leqno,draft]{article}
\usepackage{amsfonts}
\usepackage{amsmath, amsthm, amsfonts, amssymb, color}
\usepackage{mathrsfs}
\usepackage{color}
\theoremstyle{definition}

\newcommand{\scr}[1]{\mathscr #1}
\definecolor{wco}{rgb}{0.5,0.2,0.3}

\numberwithin{equation}{section} \theoremstyle{remark}

\newcommand{\ua}{\uparrow}

\title{{\bf Order Preservation for  Path-Distribution Dependent SDEs}\footnote{Supported in
 part by  NNSFC (11771326, 11431014).} }
\author{
{\bf   Xing Huang $^{a)}$,  Chang Liu $^{c)}$ Feng-Yu Wang $^{a), b)}$  }\\
\footnotesize{ a)Center of Applied Mathematics, Tianjin
University, Tianjin 300072, China}\\
 \footnotesize{ b)Department of Mathematics,
Swansea University, Singleton Park, SA2 8PP, United Kingdom}\\
\footnotesize{ c)School of Mathematical Sciences,
Beijing Normal
University, Beijing 100875, China}\\
\footnotesize{  wangfy@tju.edu.cn}}
\begin{document}
\allowdisplaybreaks
\def\R{\mathbb R}  \def\ff{\frac} \def\ss{\sqrt} \def\B{\mathbf
B} \def\W{\mathbb W}
\def\N{\mathbb N} \def\kk{\kappa} \def\m{{\bf m}}
\def\ee{\varepsilon}\def\ddd{D^*}
\def\dd{\delta} \def\DD{\Delta} \def\vv{\varepsilon} \def\rr{\rho}
\def\<{\langle} \def\>{\rangle} \def\GG{\Gamma} \def\gg{\gamma}
  \def\nn{\nabla} \def\pp{\partial} \def\E{\mathbb E}
\def\d{\text{\rm{d}}} \def\bb{\beta} \def\aa{\alpha} \def\D{\scr D}
  \def\si{\sigma} \def\ess{\text{\rm{ess}}}
\def\beg{\begin} \def\beq{\begin{equation}}  \def\F{\scr F}
\def\Ric{\text{\rm{Ric}}} \def\Hess{\text{\rm{Hess}}}
\def\e{\text{\rm{e}}} \def\ua{\underline a} \def\OO{\Omega}  \def\oo{\omega}
 \def\tt{\tilde} \def\Ric{\text{\rm{Ric}}}
\def\cut{\text{\rm{cut}}} \def\P{\mathbb P} \def\ifn{I_n(f^{\bigotimes n})}
\def\C{\scr C}      \def\aaa{\mathbf{r}}     \def\r{r}
\def\gap{\text{\rm{gap}}} \def\prr{\pi_{{\bf m},\varrho}}  \def\r{\mathbf r}
\def\Z{\mathbb Z} \def\vrr{\varrho} \def\ll{\lambda}
\def\L{\scr L}\def\Tt{\tt} \def\TT{\tt}\def\II{\mathbb I}
\def\i{{\rm in}}\def\Sect{{\rm Sect}}  \def\H{\mathbb H}
\def\M{\scr M}\def\Q{\mathbb Q} \def\texto{\text{o}} \def\LL{\Lambda}
\def\Rank{{\rm Rank}} \def\B{\scr B} \def\i{{\rm i}} \def\HR{\hat{\R}^d}
\def\to{\rightarrow}\def\l{\ell}\def\iint{\int}
\def\EE{\scr E}\def\Cut{{\rm Cut}}
\def\A{\scr A} \def\Lip{{\rm Lip}}
\def\BB{\scr B}\def\Ent{{\rm Ent}}\def\L{\scr L}
\def\R{\mathbb R}  \def\ff{\frac} \def\ss{\sqrt} \def\B{\mathbf
B}
\def\N{\mathbb N} \def\kk{\kappa} \def\m{{\bf m}}
\def\dd{\delta} \def\DD{\Delta} \def\vv{\varepsilon} \def\rr{\rho}
\def\<{\langle} \def\>{\rangle} \def\GG{\Gamma} \def\gg{\gamma}
  \def\nn{\nabla} \def\pp{\partial} \def\E{\mathbb E}
\def\d{\text{\rm{d}}} \def\bb{\beta} \def\aa{\alpha} \def\D{\scr D}
  \def\si{\sigma} \def\ess{\text{\rm{ess}}}
\def\beg{\begin} \def\beq{\begin{equation}}  \def\F{\scr F}
\def\Ric{\text{\rm{Ric}}} \def\Hess{\text{\rm{Hess}}}
\def\e{\text{\rm{e}}} \def\ua{\underline a} \def\OO{\Omega}  \def\oo{\omega}
 \def\tt{\tilde} \def\Ric{\text{\rm{Ric}}}
\def\cut{\text{\rm{cut}}} \def\P{\mathbb P} \def\ifn{I_n(f^{\bigotimes n})}
\def\C{\scr C}      \def\aaa{\mathbf{r}}     \def\r{r}
\def\gap{\text{\rm{gap}}} \def\prr{\pi_{{\bf m},\varrho}}  \def\r{\mathbf r}
\def\Z{\mathbb Z} \def\vrr{\varrho} \def\ll{\lambda}
\def\L{\scr L}\def\Tt{\tt} \def\TT{\tt}\def\II{\mathbb I}
\def\i{{\rm in}}\def\Sect{{\rm Sect}}  \def\H{\mathbb H}
\def\M{\scr M}\def\Q{\mathbb Q} \def\texto{\text{o}} \def\LL{\Lambda}
\def\Rank{{\rm Rank}} \def\B{\scr B} \def\i{{\rm i}} \def\HR{\hat{\R}^d}
\def\to{\rightarrow}\def\l{\ell}
\def\8{\infty}\def\I{1}\def\U{\scr U}
\maketitle

\begin{abstract} Sufficient and necessary conditions  are presented for the order preservation of  path-distribution dependent SDEs. Differently from  the corresponding study of distribution independent SDEs, to   investigate the necessity of  order preservation for the present model we need to construct a family of probability spaces  in terms of the   ordered pair of initial distributions.
\end{abstract} \noindent
 \vskip 2cm

\section{Introduction}

The order preservation of stochastic processes is a  crucial property for one to compare a complicated process with simpler ones, and a result to ensure this property is called $``$comparison theorem" in the literature. There are two different type order preservations, one is in the distribution (weak) sense and the other is in the pathwise (strong) sense, where the latter implies the former.
  The weak order preservation has been investigated for  diffusion-jump Markov processes in \cite{CW,W} and references within,  as well as a class of super processes in \cite{W1}. There are also plentiful results on the strong order preservation,   see, for instance,
   \cite{BY, GD, IW, M, O, PY, PZ, YMY, Z} and references within for comparison theorems on forward/backward SDEs (stochastic differential equations), with jumps and/or with  memory. Recently,   sufficient and necessary conditions have been derived in   \cite{HW} for the order preservation of  SDEs with memory.

On the other hand, to characterize  the non-linear Fokker-Planck equations, path-distribution dependent SDEs  have been investigated in \cite{HRW}, see also \cite{FYW} and references within for distribution-dependent SDEs without memory. The aim of this paper is to present sufficient and necessary conditions of the order preservations for path-distribution dependent SDEs.

Let $r_0\ge 0$ be a constant and   $d\ge 1$ be a natural number. The path space
$\C=C([-r_0,0];\R^d)$ is Polish under the uniform norm $\|\cdot\|_\infty$.
For any continuous map $f: [-r_0,\infty)\to \R^d$  and
$t\ge 0$,  let  $f_t\in\C$ be such that $f_t(\theta)=f(\theta+t)$ for $\theta\in
[-r_0,0]$. We call $(f_t)_{t\ge 0}$ the segment of $(f(t))_{t\ge -r_0}.$
Next, let $\scr P(\C)$ be the set of probability measures on $\C$ equipped with the weak topology.
Finally, let $W(t)$ be an $m$-dimensional Brownian motion on a complete
filtration probability space $(\OO, \{\F_{t}\}_{t\ge 0},\P)$.

We consider the following Distribution-dependent SDEs with memory:
\beq\label{E1} \beg{cases} \d X(t)= b(t,X_t,\L_{X_t})\,\d t+ \sigma(t,X_t,\L_{X_t})\,\d W(t),\\
\d \bar{X}(t)= \bar{b}(t,\bar{X}_t,\L_{\bar{X}_t})\,\d t+ \bar{\sigma}(t,\bar{X}_t,\L_{\bar{X}_t})\,\d W(t),\end{cases}\end{equation}
where
$$b,\bar{b}: [0,\infty)\times \C\times \scr P(\C)\to \R^d;\ \ \si,\bar{\sigma}: [0,\infty)\times \C\times \scr P(\C)\to \R^d\otimes\R^m$$
are   measurable.

For any $s\ge0$ and $\F_s$-measurable $\C$-valued random variables $\xi,\bar\xi$, a solution to (\ref{E1}) for $t\ge s$ with $(X_s,\bar X_s)= (\xi,\bar\xi)$ is a  continuous adapted process $(X(t),\bar X(t))_{t\ge s}$  such that for all $t\ge s,$
\beg{equation*}\beg{split} &X(t) = \xi(0)+ \int_s^t  b(r,X_r,\L_{X_r})\d r+ \int_s^t \si(r, X_r,\L_{X_r})\d W(r ),\\
&\bar X(t) = \bar\xi(0)+ \int_s^t \bar b(r,\bar X_r,\L_{\bar{X}_r})\d r+ \int_s^t\bar \si(r,\bar X_r,\L_{\bar{X}_r})\d W(r) ,\end{split}\end{equation*} where  $(X_t, \bar X_t)_{t\ge s}$ is the segment process of $(X(t), \bar X(t))_{t\ge s-r_0}$ with
$(X_s,\bar X_s)= (\xi,\bar\xi)$.

Following the line of  \cite{HRW}, we consider the class of probability measures of finite second moment:
$$\scr P_2(\C)=\bigg\{\nu\in \scr P(\C):\nu(\|\cdot\|_\infty^2):=\int_{\C}\|\xi\|_\infty^2\nu(d\xi)<\infty\bigg\}.$$
It is a Polish space under the Wasserstein distance
$$\mathbb{W}_2(\mu_1,\mu_2):=\inf_{\pi\in \mathbf{C}(\mu_1,\mu_2)}\left(\int_{\C \times \C}\|\xi-\eta\|_\infty^2\pi(d\xi,d\eta)\right)^{\ff 1 2},\ \ \mu_1,\mu_2\in\scr P_2(\C),$$
where $\mathbf{C}(\mu_1,\mu_2)$ is the set of all couplings for $\mu_1$ and $\mu_2;$ that is, $\pi\in \mathbf{C}(\mu_1,\mu_2)$ if it is a probability measure on $\C^2$ such that
$$\pi(\cdot\times \C)= \mu_1,\ \ \pi(\C\times\cdot)=\mu_2.$$

To investigate  the order preservation, we make   the following assumptions.

\beg{enumerate} \item[{\bf (H1)}] (Continuity)
There exists an increasing function $\aa: \R_+\to \R_+$     such that for any $ t\ge 0; \xi,\eta\in \C; \mu,\nu\in
\scr P_2(\C)$,
 \begin{align*}
&|b(t,\xi,\mu)- b(t,\eta,\nu)|^2+|\bar{b}(t,\xi,\mu)- \bar{b}(t,\eta,\nu)|^2+\|\si(t,\xi,\mu)- \si(t,\eta,\nu)\|_{HS}^2\\
&+\|\bar{\si}(t,\xi,\mu)- \bar{\si}(t,\eta,\nu)\|_{HS}^2\le \aa(t)\big(\|\xi-\eta\|_{\infty}^2+   \W_2(\mu,\nu)^2\big).
\end{align*}
\item[{\bf (H2)}] (Growth)  There exists an increasing function $K: \R_+\to\R_+   $ such that
 \begin{align*}
  |b(t,0,\dd_0)|^2+ |\bar b(t,0,\dd_0)|^2+\|\si(t,0,\dd_0)\|_{HS}^{2}+\|\bar\si(t,0,\dd_0)\|^{2}_{HS}
 \le K (t),\ \ t\ge 0,
 \end{align*}where $\dd_0$ is the Dirac measure at point $0\in\C$.
\end{enumerate}

It is easy to see that these two conditions imply assumptions $(H1)$-$(H3)$ in   \cite{HRW},  so that by \cite[Theorem 3.1]{HRW},  for any $s\ge 0$ and $\F_s$-measurable $\C$-valued random variables $\xi,\bar\xi$ with finite second moment,
 the equation  (\ref{E1}) has a unique solution  $\{X(s,\xi;t), \bar X(s,\bar \xi;t)\}_{t\ge s}$    with $X_s=\xi$ and $\bar X_s=\bar\xi$.
Moreover, the segment process $\{X(s,\xi)_t, \bar X(s,\bar\xi)_t\}_{t\ge s}$ satisfies
\beq\label{BDD} \E \sup_{t\in [s,T]} \big(\|X(s,\xi)_t\|_\infty^2+ \|\bar X(s,\bar \xi)_t\|_\infty^2\big)<\infty,\ \ T\in [s,\infty).\end{equation}

To characterize the order-preservation for solutions of \eqref{E1}, we introduce the partial-order on $\C.$ For $x=(x^1,\cdots, x^d)$ and $ y=(y^1,\cdots, y^d)\in\R^d$, we write $x\le y$ if $x^i\le y^i$ holds for all $1\le i\le d.$
Similarly,  for $\xi=(\xi^1,\cdots,\xi^d)$ and $\eta=(\eta^1,\cdots,\eta^d)\in\C$, we write $\xi\le \eta$ if $\xi^i(\theta)\le \eta^i(\theta)$ holds for all $\theta\in [-r_0,0]$ and $1\le i\le d.$ A function $f$ on $\C$ is called increasing if $f(\xi)\le f(\eta)$ for $\xi\le \eta$. Moreover, for any $\xi_1,\xi_2\in\C$, $\xi_1\land\xi_2\in \C$ is defined by $$(\xi_1\land\xi_2)^i:=\min\{\xi^i_1,\xi_2^i\}, \ \  1\le i\le d.$$
For two probability measures $\mu,\nu\in\scr P(\C)$, we   write $\mu\leq \nu$    if $\mu(f)\leq \nu(f)$  holds for any increasing function $f\in C_b(\C)$.  According to \cite[Theorem 5]{KKO},    $\mu\leq \nu$  if and only if there exists $\pi\in \mathbf C(\mu,\nu)$ such that
 $\pi(\{(\xi,\eta)\in \C^2: \xi\le \eta\})=1.$

\beg{defn} The stochastic differential system  $(\ref{E1})$ is called order-preserving, if  for any $s\ge 0$ and $\xi, \bar\xi\in L^2(\OO\to\C, \F_{s},\P)$ with $\P(\xi\le\bar\xi)=1$,   $$\P\big(X(s,\xi;t)\le \bar X(s,\bar\xi;t),\ t\ge s\big)=1.$$   \end{defn}

We first present the following sufficient conditions for the order preservation, which reduce back to the corresponding ones in \cite{HW} when the system is distribution-independent.

\beg{thm}\label{T1.1} Assume {\bf (H1)} and {\bf (H2)}. The system \eqref{E1} is order-preserving provided the following  two conditions are satisfied:
\beg{enumerate}
\item[$(1)$] For any  $1\leq i\leq d$, $\mu,\nu\in \scr P_2(\C)$ with $\mu\leq\nu$, $\xi,\eta\in\C$ with $\xi\leq \eta$ and $\xi^i(0)=\eta^i(0)$,
    \begin{equation*}
b^i(t,\xi,\mu)\leq \bar{b}^i(t, \eta, \nu),\ \ \text{a.e.}\ t\ge 0.
    \end{equation*}
\item[$(2)$] For a.e.\ $t\ge 0$ it holds: $\si(t, \cdot,\cdot)= \bar\si(t,\cdot,\cdot)$    and $\sigma^{ij}(t,\xi,\mu)=\si^{ij}(t, \eta,\nu)$  for any $1\leq i\leq d$, $1\leq j\leq m$,  $\mu,\nu\in \scr P_2(\C)$ and $\xi,\eta\in \C$ with $\xi^{i}(0)=\eta^i(0)$.
\end{enumerate}
\end{thm}

Condition (2) means that for a.e. $t\ge 0$, $\si(t,\xi,\mu)=\bar\si(t,\xi,\mu)$   and the  dependence of $\si^{ij}(t,\xi,\mu)$ on $(\xi,\mu)$ is only via $\xi^i(0)$.

On the other hand, the next result shows that    these conditions are also necessary if all coefficients are continuous on $[0,\infty)\times\C\times\scr P_2(\C)$, so that \cite[Theorem 1.2]{HW} is covered when the system is distribution-independent.

\beg{thm}\label{T1.2} Assume {\bf (H1)},  {\bf (H2)} and that $\eqref{E1}$ is order-preserving  for any complete filtered probability space $(\OO,\{\F_t\}_{t\ge 0},\P)$ and $m$-dimensional Brownian motion $W(t)$ thereon.  Then  for any $1\leq i\leq d$, $\mu,\nu\in \scr P_2(\C)$ with $\mu\le \nu$, and $\xi,\eta\in\C$ with $\xi\le \eta$ and $\xi^i(0)=\eta^i(0)$, the following assertions hold:
\beg{enumerate}
\item[$(1')$]  $b^i(t,\xi,\mu)\leq \bar{b}^i(t, \eta, \nu)$ if  $b^i$ and $\bar b^i$ are continuous  at points $(t,\xi,\mu)$ and $(t, \eta,\nu)$ respectively.
\item[$(2')$] For any $1\le j\le m$, $\sigma^{ij}(t,\xi,\mu)=\bar\si^{ij}(t, \eta,\nu)$ if $\si^{ij}$ and $\bar \si^{ij}$ are continuous   at points $(t,\xi,\mu)$ and $(t, \eta,\nu)$ respectively.
\end{enumerate}
Consequently, when $b, \bar b, \si$ and $\bar \si$ are continuous on $[0,\infty)\times \C\times \scr P_2(\C)$, conditions $(1)$ and $(2)$ hold.
\end{thm}

These two theorems will be proved in Section 2 and Section 3 respectively.

\section{Proof  of Theorem  \ref{T1.1}}

  Assume  {\bf (H1)} and {\bf (H2)}, and let conditions (1) and (2) hold. For any $T>t_0\ge 0$ and $\xi, \bar\xi\in L^2(\OO\to\C, \F_{t_0},\P)$ with $\P(\xi\le\bar\xi)=1$, it suffices  to prove
  \beq\label{P} \E\sup_{t\in[t_0,T]} (X^i(t_0,\xi; t)-\bar X^i(t_0,\bar\xi;t))^+=0,\ \ 1\le i\le d,\end{equation} where $s^+:=\max\{0,s\}.$
  For simplicity, in the following we   denote $X(t)=X(t_0,\xi;t)$ and $\bar X(t)=\bar X(t_0,\bar\xi;t)$ for $t\ge t_0-r_0$. Then
  $$X(t)=\xi(t-t_0),\ \ \bar X(t)= \bar \xi(t-t_0),\ \ t\in [t_0-r_0,t_0].$$

To prove   \eqref{P} using It\^o's formula, we    take  the following $C^2$-approximation of $s^+$  as in the proof of \cite[ Theorem 1.1]{HW}.
For any $n\ge 1$, let $\psi_n: \R\to [0,\infty)$ be constructed as follows: $\psi_n(s)=\psi_n'(s)=0$ for $s\in (-\infty,0]$, and
$$\psi_n''(s)=\beg{cases} 4n^2s, & s\in [0,\ff 1 {2n}],\\
-4n^2(s-\ff 1 n), & s\in [\ff 1 {2n}, \ff 1 n],\\
0, &\text{otherwise}.\end{cases}$$  We have
\beq\label{1.3} 0\le \psi_n'\le 1_{(0,\infty)}, \ \text{and\ as\ } n\uparrow\infty: \ 0\le \psi_n(s)\uparrow s^+,\ \ s\psi_n''(s)\le  1_{(0,\ff 1 n)}(s)\downarrow 0.\end{equation}
Let $$\tau_k=\inf\big\{t\ge t_0: |X(t)-X(t)\land\bar X(t)|\ge k\big\},\ \ k\ge 1.$$ Since $$\psi_n(X^i(t_0)-\bar X^i(t_0))=\psi_n(\xi^i(0)-\bar\xi^i(0))=0,$$  and due to (2)   $\si(t,\cdot,\cdot)=\bar\si(t,\cdot,\cdot)$ for a.e. $t\ge 0$, by
 It\^o's  formula we obtain
\beq\label{1.4} \beg{split}&\psi_n(X^i(t\land \tau_k)-\bar X^i(t\land\tau_k))^2\\
&= M_i(t\land\tau_k)+2\int_{t_0}^{t\land\tau_k} (b^i(s,X_s,\L_{X_s})-\bar b^i(s,\bar X_s,\L_{\bar{X}_s}))\{\psi_n\psi_n'\}(X^i(s)-\bar X^i(s))\d s\\
&\quad +  \sum_{j=1}^m \int_{t_0}^{t\land\tau_k} (\si^{ij}(s,X_s,\L_{X_s})- \si^{ij}(s,\bar X_s,\L_{\bar{X}_s}))^2\{ \psi_n\psi_n''+\psi_n'^2\}(X^i(s)-\bar X^i(s))\d s \end{split}\end{equation}for any $k,n\ge 1$, $1\le i\le d$ and $t\ge t_0,$
where
$$M_i(t):= 2\sum_{j=1}^m \int_{t_0}^t (\si^{ij}(s, X_s,\L_{X_s})- \si^{ij}(s,\bar X_s, \L_{\bar{X}_s}))\{\psi_n\psi_n'\}(X^i(s)-\bar X^i(s))\d W^j(s).$$
Noting that $0\le\psi_n'(X^i(s)-\bar X^i(s))\le 1_{\{X^i(s)>\bar X^i(s)\}}$ and when $X^i(s)>\bar X^i(s)$ one has
$(X_s\land \bar X_s)^i(0)=(\bar X_s)^i(0),$ it follows from (1) that  for a.e. $s\in [t_0,T],$
\begin{equation*}\begin{split}
&[b^i(s,X_s\land\bar X_s,\L_{X_s\land\bar X_s})-\bar b^i(s,\bar X_s,\L_{\bar X_s})]{\psi_n\psi_n'}(X^i(s)-\bar X^i(s))\leq0,\ \ n\ge 1.
\end{split}\end{equation*}
 Combining this with {\bf (H1)} and $0\le\psi_n'\le 1$, we obtain
\beq\label{1.5} \beg{split}
&2\int_{t_0}^{t\land\tau_k}[b^i(s,X_s,\L_{X_s})-\bar b^i(s,\bar X_s,\L_{\bar{X}_s})]\{\psi_n\psi_n'\}(X^i(s)-\bar X^i(s))\d s\\
& =2\int_{t_0}^{t\land\tau_k}\Big\{[b^i(s,X_s,\L_{X_s})-b^i(s,X_s\land \bar X_s,\L_{X_s\land \bar X_s})]\{\psi_n\psi_n'\}(X^i(s)-\bar X^i(s))\\
&\quad+ [b^i(s,X_s\land\bar X_s,\L_{X_s\land \bar X_s})-\bar b^i(s,\bar X_s, \L_{\bar X_s})]\{\psi_n\psi_n'\}(X^i(s)-\bar X^i(s))\Big\}\d s\\
&\leq2\int_{t_0}^{t\land\tau_k}[b^i(s,X_s,\L_{X_s})- b^i(s,X_s\land \bar X_s,\L_{X_s\land \bar X_s})]\{\psi_n\psi_n'\}(X^i(s)-\bar X^i(s))\d s\\
&\leq\int_{t_0}^{t\land\tau_k}[|b^i(s,X_s,\L_{X_s})- b^i(s,X_s\land \bar X_s,\L_{X_s\land \bar X_s})|^2+\psi_n(X^i(s)-\bar X^i(s))^2]\d s\\
&\leq\int_{t_0}^{t\land\tau_k}\aa(T)\left[\|X_s-X_s\land \bar X_s\|_{\infty}^2+ \W_2(\L_{X_{s\wedge \tau_k}},\L_{X_{s\wedge \tau_k}\land \bar X_{s\wedge \tau_k}})^2\right]\d s\\
&\quad+\int_{t_0}^{t\land\tau_k}\psi_n(X^i(s)-\bar X^i(s))^2\d s,\ \ n\geq1,t\in[t_0,T].
\end{split}\end{equation}
Next, by (2), for a.e. $s\in [t_0,T]$, $\si^{ij}(s, X_s,\L_{X_s})=\bar\si^{ij}(s, X_s,\L_{X_s})$ depends only on $X^i(s)$. So,   \eqref{1.3} and {\bf (H1)} yield
\beq\label{1.6} \beg{split}&\sum_{j=1}^m\int_{t_0}^{t\land \tau_k}  |\si^{ij}(s,X_s,\L_{X_s})- \si^{ij}(s,\bar X_s,\L_{\bar{X}_s})|^2\{\psi_n\psi_n''+\psi_n'^2\}(X^i(s)-\bar X^i(s))\d s \\
&\leq  \int_{t_0}^{t\land \tau_k} (1_{\{X^i(s)-\bar X^i(s)\in (0,\ff 1 n)\}}+1) \aa(T)\{(X^i(s)- \bar X^i(s))^+\}^2\d s  \\
&\leq2\aa(T)\int_{t_0}^{t\land \tau_k} \{(X^i(s)- \bar X^i(s))^+\}^2\d s,\ \ n\ge 1, t\in [t_0,T],
\end{split}\end{equation}
and
\beq\label{1.7}  \beg{split} &\sum_{j=1}^m\int_{t_0}^{t\land \tau_k}   \big|\si^{ij}(s, X_s,\L_{X_s})- \si^{ij}(s,\bar X_s, \L_{\bar{X}_s})\big|^2\big|(\psi_n\psi_n')(X^i(s)-\bar X^i(s))\big|^2\d s \\
&\le \alpha(T)\int_{t_0}^{t\land \tau_k}  \{(X^i(s)- \bar X^i(s))^+\}^2  \psi_n(X^i(s)-\bar X^i(s))^2\d s,\ \ n\ge 1, t\in [t_0,T].\end{split}\end{equation} By the Burkholder-Davis-Gundy  inequality, we deduce from \eqref{1.7} that
\beq\label{QQ} \beg{split}  \E\sup_{s\in [t_0,t]}  M_i(s\land \tau_k)
  \le & C_1  \E \bigg(\int_{t_0}^{t\land\tau_k}  \{(X^i(s)- \bar X^i(s))^+\}^2  \psi_n(X^i(s)-\bar X^i(s))^2\d s\bigg)^{\ff 1 2}\\
 \le &2 C_1^2  \E \int_{t_0}^{t\land\tau_k} \|X_s-X_s\land \bar X_s\|_\infty^2 \d s\\
 &+  \ff 1 8 \E\sup_{s\in [t_0, t\land\tau_k]}\psi_n(X^i(s)-\bar X^i(s))^2,\ \ n\ge 1, t\in [t_0,T]  \end{split}
\end{equation} holds for some  constant  $C_1>0.$

Now, let
$$\phi_k(s)=\sup_{r\in [t_0-r_0, s\land\tau_k]}|X(r)-X(r)\land \bar X(r)|^2,\ \ s\in [t_0,T].$$
By substituting \eqref{1.5}, \eqref{1.6} and \eqref{QQ} into \eqref{1.4} and noting that  $X_{t_0}\le\bar X_{t_0}$ and $$\W_2(\L_{X_{s\wedge \tau_k}},\L_{X_{s\wedge \tau_k}\land \bar X_{s\wedge \tau_k}})^2\leq \mathbb{E}\|X_{s\wedge \tau_k}-X_{s\wedge \tau_k}\land \bar X_{s\wedge \tau_k}\|_{\infty}^2,$$ we obtain
\beq\label{1.8} \beg{split}
&\E\sup_{r\in [t_0-r_0, t\land \tau_k]} \psi_n(X^i(r)-\bar X^i(r))^2=\E\sup_{r\in [t_0, t\land \tau_k]} \psi_n(X^i(r)-\bar X^i(r))^2\\
&\leq C_2\E\int_{t_0}^{t\land\tau_k}\|X_{s\wedge \tau_k}-X_{s\wedge \tau_k}\land \bar X_{s\wedge \tau_k}\|_{\infty}^2\d s+\ff 1 8 \E\sup_{s\in [t_0, t\land\tau_k]}\psi_n(X^i(s)-\bar X^i(s))^2,
\end{split}\end{equation}
for some constants $C_2>0$ and all $k,n\ge 1, t\in [t_0,T], 1\le i\le d.$ Therefore, there exists a constant $C>0$ such that for any $n,k\ge 1$ and $t\in [t_0,T],$
\beg{equation*}\beg{split} &  \sum_{i=1}^d\E\sup_{r\in [t_0-r_0, t\land \tau_k]} \psi_n(X^i(r)-\bar X^i(r))^2\le C\int_{t_0}^{t}\E\phi_k(s)\d s,\ \ k,n\ge 1.\end{split}\end{equation*}
 Letting  $n\uparrow\infty$, we arrive at
$$\E  \phi_k(t)\le C  \int_{t_0}^t  \E\phi_k(s)\d s,\ \ t\in [t_0,T], k\ge 1.$$ By the definition of $\phi_k$ and $\tau_k$, $\E\phi_k(t)$ is locally bounded in $t\ge 0$. So, Gronwall's inequality  implies
$$ \E \phi_k(T)=0,\ \ k\ge 1.$$ Letting $k\uparrow\infty$  we  prove  (\ref{P}).

\section{Proof of Theorem \ref{T1.2}}

 We first observe that when $b, \bar b, \si, \bar \si$ are continuous on $[0,\infty)\times \C\times\scr P_2(\C)$, $(1')$ and $(2')$ imply $(1)$ and $(2)$. Obviously, $(1')$ implies $(1)$. It suffices to prove $(2)$.

 Firstly, taking $\xi=\eta$ and $\mu=\nu$, by the continuity of $\si$ and $\bar \si$, $(2')$ implies $\si=\bar\si$.

 In general, for  $1\leq i\leq d$, $\mu,\nu\in \scr P_2(\C)$,   and $\xi,\eta\in\C$ with  $\xi^i(0)=\eta^i(0)$,
 we take
 $$(\mu\land \nu) (\cdot)= (\mu\times\nu)\big(\{(\xi_1,\xi_2)\in \C^2: \xi_1\land \xi_2\in \cdot\}\big).$$

 Then $\mu\land \nu\in \scr P_2(\C)$ and $\mu\land \nu\le \mu,\ \mu\land \nu\le\nu$. Moreover, $\xi\land\eta\le \xi,\ \xi\land\eta\le \eta$ with
    $\xi^i(0)=(\xi\land\eta)^i(0)=\eta^i(0)$. So, applying  $(2')$ twice we obtain
 $$\si^{ij}(t,\xi,\mu)=\si^{ij}(t,\xi\land \eta, \mu\land\nu)= \si^{ij}(t, \eta,\nu).$$
Since $\si=\bar\si$, this implies  $(2)$.

\

 Now,  let $t_0\ge 0$  $1\leq i\leq d$, $\mu,\nu\in \scr P_2(\C)$ with $\mu\le \nu$,   and $\xi,\eta\in\C$ with $\xi\le\eta$ and $\xi^i(0)=\eta^i(0)$.  To prove $(1')$ and $(2')$ for $t=t_0$, we   construct a family of  complete filtration probability spaces $(\OO,\{\F_t^\vv\}_{t\geq 0},\P^\vv)_{\vv\in [0,1)}$, $m$-dimensional Brownian motion $W(t)$, and initial random variables $ X_{t_0}\le \bar X_{t_0}$ as follows.

Firstly, since $\mu\le \nu$, by \cite[Theorem 5]{KKO} we may take $\pi_0\in {\mathbf C}(\mu,\nu)$ such that
\beq\label{OS1}\pi_0(\{(\xi_1,\xi_2)\in \C^2: \xi_1\le \xi_2\})=1.\end{equation} For any $\vv\in [0,1)$, let
\beq\label{OS2}\pi_\vv= (1-\vv)\pi_0+\vv \dd_{(\xi,\eta)},\end{equation}where $\dd_{(\xi,\eta)}$ is the Dirac measure at point $(\xi,\eta)$.
Let $\P_0$ be the standard Wiener measure on $\OO_0:=C([0,\infty)\to \mathbb{R}^m)$, and let $\F_{0,t}$ be the completion of $\si(\oo_0\mapsto \oo_0(s): s\le t)$ with respect to the Wiener measure. Then
the coordinate process $\{W_0(t)\}(\oo_0):= \oo_0(t),\ \oo_0\in\OO_0, t\ge 0$ is an $m$-dimensional Brwonian motion on the filtered probability space $(\OO_0, \{\F_{0,t}\}_{t\ge 0}, \P_0)$.

Next, for any $\vv\in[0,1)$, let
$\OO=\OO_0\times \C^2,\ \P^\vv= \P_0 \times \pi_\vv$ and $\F_t^\vv$ be the completion of $\F_{0,t}\times \B(\C^2)$ under the probability measure  $\P^\vv$. Then the process
$$\{W(t)\}(\oo):= \{W_0(t)\}(\oo_0)=\oo_0(t),\ \ t\ge 0, \oo=(\oo_0;\xi_1,\xi_2)\in \OO=\OO_0\times \C^2$$
is an $m$-dimensional Brownian motion on the complete probability space $(\OO, \{\F_t^\vv\}_{t\ge 0}, \P^\vv)$.

Finally, let
$$X_{t_0}(\oo):= \xi_1,\ \ \bar X_{t_0}(\oo):= \xi_2,\ \ \oo=(\oo_0; \xi_1,\xi_2)\in \OO=\OO_0\times \C^2.$$
They are $\F_{t_0}^\vv$-measurable random variables with
\beq\label{XT}  \scr L_{X_{t_0}}|_{\P^\vv}=\mu_\vv:=\pi_\vv(\cdot\times\C),\ \ \  \scr L_{\bar X_{t_0}}|_{\P^\vv}=\nu_\vv:=\pi_\vv(\C\times\cdot),\end{equation} where $\scr L_X|_{\P^\vv}$ denotes the distribution of a random variable $X$ under probability $\P^\vv$. By $\xi\le \eta$ and \eqref{OS1}, \eqref{OS2}, we have
\beq\label{XT0} \P^\vv(X_{t_0}\le \bar X_{t_0})=\pi_\vv\big(\{(\xi_1,\xi_2)\in \C^2:\ \xi_1\le \xi_2\}\big)=1,\ \ \vv\in [0,1).\end{equation}
So, letting $(X_t,\bar X_t)_{t\ge t_0}$ be the segment process of the solution to \eqref{E1} with initial value $(X_{t_0}, \bar X_{t_0})$,      the order preservation implies
\beq\label{OP}\P^\vv\big( X_t\le \bar X_t,\ t\ge t_0\big)=1,\ \ \vv\in [0,1).\end{equation}
Let $\E^\vv$ be the expectation for $\P^\vv$. With the above preparations, we are able to prove
  $(1')$ and $(2')$ as follows.

\beg{proof}[Proof of $(1')$]  Let  $b^i, \bar b^i$ be continuous at points $(t_0, \xi,\mu)$ and $(t_0,\eta,\nu)$ respectively.
We intend to prove
 $b^i(t_0, \xi,\mu)\le \bar{b}^i(t_0,\eta,\nu)$. Otherwise, there exists a constant $c_0>0$ such that
 \beq\label{CB} b^i(t_0, \xi,\mu)\ge c_0+ \bar b^i(t_0,\eta,\nu).\end{equation}
 Let $\mu_\vv,\nu_\vv$ be in \eqref{XT}.
 Obviously, $\{\mu_\vv,\nu_\vv\}_{\vv\in [0,1)}$ are bounded in $\scr P_2(\C)$ and,  as $\vv\to 0$,  $\mu_\vv\to\mu,\ \nu_\vv\to\nu$ weakly.
Consequently,
 $$\lim_{\vv\downarrow 0} \{\W_2(\mu_\vv,\mu)+\W_2(\nu_\vv,\nu)\}=0.$$
 Combining this with {\bf (H1)} and \eqref{CB}, there exists $\vv\in (0,1)$ such that
 \beq\label{CB2} b^i(t_0, \xi,\mu_\vv)\ge \ff 1 2 c_0+  \bar b^i(t_0,\eta,\nu_\varepsilon)> \bar b^i(t_0,\eta,\nu_\varepsilon).\end{equation}

 Now, consider the event
\beq\label{AX} A:=\{X_{t_0}= \xi, \bar X_{t_0}=\eta\}\in \F_{t_0}^\vv.\end{equation}
 Then
 \beq\label{XT3} \P^\vv(A)\ge \vv \dd_{(\xi,\eta)}(\{(\xi,\eta)\})=\vv>0.\end{equation}
By \eqref{E1} and \eqref{OP}, for any $s\ge 0$, $\P^\vv$-a.s.
\beq\label{XA}\beg{split}
0&\ge X^i(t_0+s)-\bar X^i(t_0+s)   =X^i_{t_0}(0)-\bar X^i_{t_0}(0)\\
 &\quad +   \int_{t_0}^{t_0+s}b^i(r,X_r,\L_{X_r}|_{\P^\vv})\,\d r-\int_{t_0}^{t_0+s}\bar b^i(r,\bar X_r,\L_{\bar{X}_r}|_{\P^\vv})\,\d r  \\
&\quad +  \sum_{j=1}^m \int_{t_0}^{t_0+s}\sigma^{ij}(r,X_r,\L_{X_r}|_{\P^\vv})\,\d W^j(r)-\int_{t_0}^{t_0+s}\bar \sigma^{ij}(r,\bar X_r,\L_{\bar{X}_r}|_{\P^\vv})\,\d W^j(r).\end{split}\end{equation}
 Since  on $A$ we have $X^i_{t_0}(0)-\bar X^i_{t_0}(0)=0$, and $A$ is $\F_{t_0}^\vv$-measurable, by multiplying \eqref{XA} with $\ff {1_A} s$ and  taking conditional expectation  with respect to $\F^\vv_{t_0}$,   we
obtain    $\P^\vv$-a.s.
$$\E^\vv\bigg(\ff{1_A}s \int_{t_0}^{t_0+s}b^i(r,X_r,\L_{X_r}|_{\P^\vv})\,\d r\bigg|\F^\vv_{t_0}\bigg) \leq \E^\vv\bigg(  \ff{1_A}s\int_{t_0}^{t_0+s}\bar{b}^i(r,\bar{X}_r,\L_{\bar{X}_r}|_{\P^\vv})\,\d r\bigg|\F^\vv_{t_0}\bigg),\ \ s>0.$$
Combining this with the fact that  $b^i$ and $\bar b^i$ are continuous at points   $(t_0,\xi,\mu)$ and $(t_0, \eta,\nu)$ respectively,
 and using  {\bf(H1)}, {\bf(H2)}, \eqref{BDD} and the continuity of the solution, taking $s\downarrow 0$ we derive  $\P^\vv$-a.s.
  $$
\E^\vv\big(1_Ab^i(t_0,X_{t_0}, \scr L_{X_{t_0}}|_{\P^\vv}) \big|\F^\vv_{t_0}\big) \leq\E^\vv\big(1_A\bar b^i(t_0,\bar X_{t_0}, \scr L_{\bar X_{t_0}}|_{\P^\vv}) \big|\F^\vv_{t_0}\big).$$
This together with \eqref{AX} and \eqref{XT} leads to  $\P^\vv$-a.s.
\begin{equation*}\begin{split}
b^i(t_0,\xi,\mu_\vv)1_{A}\leq \bar{b}^i(t_0,\eta,\nu_\vv)1_A,
\end{split}\end{equation*} which is impossible  according to \eqref{CB2} and \eqref{XT3}.
Therefore, $b^i(t_0, \xi,\mu)\le \bar b^i(t_0,\eta,\nu)$ has to be true.
\end{proof}

\beg{proof}[Proof of $(2')$]  Let $\si^{ij}$ and $\bar \si^{ij}$ be continuous at points $(t_0, \xi,\mu)$ and $(t_0,\eta,\nu)$ respectively.
If $\si^{ij}(t_0, \xi,\mu)\ne \bar\si^{ij}(t_0,\eta,\nu)$,   by {\bf (H1)}, there exist constants $c_1>0$ and $\vv\in (0,1)$ such that
\beq\label{SXT} |\si^{ij}(t_0, \xi, \mu_\vv)-\bar\si^{ij}(t_0, \eta,\nu_\vv)|^2\ge 2c_1>0.\end{equation}
Let
\beg{align*}&\tau=\inf\big\{t\ge t_0: |\si^{ij}(t_0, \xi, \mu_\vv)-\bar\si^{ij}(t_0, \eta,\nu_\vv)|^2\le c_1\big\},\\
&\tau_n= \tau\land\inf\big\{ t\ge t_0: |b^i(t, X_t, \scr L_{X_t}|_{\P^\vv})- \bar b^i(t, \bar X_t, \scr L_{\bar X_t}|_{\P^\vv})|^2\ge n\big\},\ \ n\ge 1.\end{align*} Let $g_n(s)=\e^{ns}-1.$ Then $g_n\in C_b^2((-\infty,0])$. By the order preservation we have $X^i_t\le \bar X^i_t,\ t\ge t_0$. So, applying It\^o's formula, we obtain $\P^\vv$-a.s.
  \beg{align*}  &0\ge  \E^\vv  \big(g_n((X^i -\bar X^i)( t\land {\tau}_n))|\F^\vv_{t_0}\big) =
  g_n((X^i -\bar X^i)(t_0))\\
&+\E^\vv \bigg(\sum_{j=1}^m\int_{t_0}^{t\land{\tau}_n} g_n'((X^i-\bar X^i)(s))\big(\sigma^{ij}(s,X_s,\L_{X_s}|_{\P^\vv})-\bar \sigma^{ij}(s,\bar X_s,\L_{\bar{X}_s}|_{\P^\vv})\big)\d W^j(s)\bigg|\F^\vv_{t_0}\bigg) \\
& + \E^\vv \bigg(\int_{t_{0}}^{t\land{\tau}_n} \bigg\{g_n'((X^i-\bar X^i)(s)) \big(b^i(s,X_s,\L_{X_s}|_{\P^\vv})-\bar b^i(s,\bar X_s,\L_{\bar{X}_s}|_{\P^\vv})\big)\\
&\ +\ff {g_n''(X^i(s)-\bar X^i(s))} 2\sum_{j=1}^m \big|\si^{ij}(s, X_s,\L_{X_s}|_{\P^\vv})-\bar\si^{ij}(s,\bar X_s,\L_{\bar{X}_s}|_{\P^\vv})\big|^2\bigg\}\d s\bigg|\F^\vv_{t_0}\bigg) \\
&\ge g_n((X^i -\bar X^i)(t_0)) +\Big(\ff{n^2c_1}{2} -  n^{3/2}\Big)\E^\vv (t\land{\tau}_n-t_0|\F^\vv_{t_0}),\ \ n\ge 1.\end{align*}
By \eqref{AX} and $\xi^i(0)=\eta^i(0),$ this implies
$$ 1_A \Big(\ff{n^2c_1}{2} -  n^{3/2}\Big)\E^\vv (t\land{\tau}_n-t_0|\F^\vv_{t_0})\le -1_A g_n((X^i-\bar X^i)(t_0) ) =-1_A g_n(0)=0$$ for all $n\ge 1$ and $t>t_0$. Therefore,
\beq\label{LK}1_A \E^\vv (t\land{\tau}_n-t_0|\F_{t_0}^\vv)= 0,\ \ \ss n> \ff{2}{c_1}, t>t_0.\end{equation}  But by ${\bf (H1)}$, \eqref{SXT}  and the continuity of the solution, on the set $A$ we have
$$\lim_{n\to\infty} \tau_n= \tau >t_0.$$ So, \eqref{LK} implies $\P^\vv(A)=0$, which contradicts \eqref{XT3}. Hence,  $\si^{ij}(t_0, \xi,\mu)= \bar \si^{ij}(t_0, \eta,\nu).$
\end{proof}


\end{document}